\documentstyle[12pt]{article}
\title{ The Radicals of Hopf Module Algebras
 \thanks {This work was supported by the  National Natural Science Foundation  }}
\author{
Shouchuan Zhang  \\ Department  of Mathematics, Hunan University\\
Changsha  410082, \
 P.R.China. \ \
E-mail:z9491@yahoo.com.cn\\
}
\date{}
\begin{document}
\newtheorem{Theorem}{\quad Theorem}[section]
\newtheorem{Proposition}[Theorem]{\quad Proposition}
\newtheorem{Definition}[Theorem]{\quad Definition}
\newtheorem{Corollary}[Theorem]{\quad Corollary}
\newtheorem{Lemma}[Theorem]{\quad Lemma}
\newtheorem{Example}[Theorem]{\quad Example}
\maketitle \addtocounter{section}{-1}

 \begin {abstract} The characterization of $H$-prime radical is given
in many ways. Meantime, the relations between the radical of smash
product $R \# H$ and the $H$-radical of Hopf module algebra $R$
are obtained.

 \end {abstract}

\section {Introduction and Preliminaries}

In this paper, let $k$ be a commutative associative ring with
unit, $H$ be an algebra with unit and comultiplication
$\bigtriangleup$ ( i.e. $\Delta $ is a linear map: $H \rightarrow
H \otimes H$),
  $R$ be an algebra
over $k$ ($R$ may be without unit) and $R$ be an $H$-module
algebra.

We define some necessary concept as follows.

If   there exists a linear map $
\left \{ \begin {array} {ll} H \otimes R & \longrightarrow R \\
 h \otimes r & \mapsto   h \cdot r \end {array} \right. $ such that
     $$h\cdot rs = \sum (h_1\cdot r)(h_2 \cdot s) \hbox { \ \ and \ }
     1_H \cdot r = r $$
     for all $r, s \in R, h\in H,$  then we say that
     $H$ weakly acts on $R.$
  For any ideal $I$ of $R$, set
  $$(I:H):= \{ x \in R \mid h\cdot x \in I \hbox { for all } h \in H  \}. $$
 $I$ is called an $H$-ideal if $h\cdot I
\subseteq I$ for any $h \in H$. Let $I_H$ denote the maximal
$H$-ideal of $R$ in $I$. It is clear that $I_H= (I:H).$ An
$H$-module algebra $R$ is called an $H$-simple module algebra if
$R$ has not any non-trivial $H$-ideals and $R^2 \not=0.$
 $R$ is said to be  $H$-semiprime if there are no
non-zero nilpotent $H$-ideals in $R$.  $R$ is said to be
  $H$-prime
if $IJ =0$ implies $I=0$ or $J=0$ for any $H$-ideals $I$ and $J$
of $R$. An $H$-ideal $I$ is called an   $H$-(semi)prime ideal of
$R$ if $R/I$ is  $H$-(semi)prime.
   $ \{ a_{n} \}$  is called an
  $H$-$m$-sequence  in  $R$ with beginning $a$
  if there exist $h_n, h_n' \in H$ such that $a_1 = a \in R$ and
  $a_{n+1} =
  (h_{n}.a_{n})b_{n}(h_{n}'.a_{n})$   for any natural number $n$.
If every $H$-$m$-sequence $\{ a_{n} \}$
  with  $a_{7.1.1} = a$, there exists a natural number $k$ such that $a_{k}= 0,$
  then $a$ is called an $H$-$m$-nilpotent element.
  Set
  $$W_{H}(R) = \{ a \in R \mid a \hbox { \ is an \ } H \hbox{-}m
  \hbox {-nilpotent element} \}. $$
 $R$ is called an $H$-module algebra
 if the following
conditions hold:

   (i)  $R$ is a unital left $H$-module(i.e. $R$ is a left $H$-module and
   $1_H \cdot a = a$
   for any $a \in R$);

   (ii)  $h\cdot ab = \sum (h_1\cdot a)(h_2 \cdot b)$ for any $a, b \in R$,
      $h\in H$, where
$\Delta (h) = \sum h_1 \otimes h_2$.  \\
$H$-module algebra is sometimes called a Hopf module.

If $R$ is an $H$-module algebra with a unit $1_R$, then
 $$h \cdot 1_R = \sum_h (h_1 \cdot 1_R)(h_2S(h_3)\cdot 1_R)$$
 $$ = \sum _h h_1 \cdot
 (1_R (S(h_2)\cdot 1_R)) = \sum _h h_1S(h_2) \cdot 1_R= \epsilon (h)1_R, $$
 i.e.         \ \ \ \ $ \hbox { \ \ \ \ \ }   h \cdot 1_R = \epsilon (h)1_R$ \\
for any $h \in H.$

An $H$-module algebra $R$ is called a unital $H$-module algebra if
$R$ has a unit $1_R$ such that $h \cdot 1_R = \epsilon (h)1_R$ for
any $h \in H$. Therefore, every $H$-module algebra with unit is a
unital $H$-module algebra. A left $R$-module $M$ is called an
$R$-$H$-module if $M$ is also a left unital $H$-module
 with $h  (am)= \sum (h_1 \cdot a)(h_2m)$  for all
$h \in H, a \in R, m \in M$. An $R$-$H$-module $M$ is called an
$R$-$H$- irreducible module if there are no non-trivial
$R$-$H$-submodules in $M$ and $RM \not=0$. An algebra homomorphism
$\psi: R \rightarrow R'$ is called an $H$-homomorphism  if $\psi
(h \cdot a) = h \cdot \psi (a)$ for any $h \in H, a \in R.$ Let
$r_b, r_j, r_l, r_{bm}$  denote the Baer radical, the Jacobson
radical, the locally nilpotent radical,
 the Brown-MacCoy radical
  of algebras respectively.   Let
   $I \lhd_H R$ denote that $I$ is an $H$-ideal of $R.$

\section{ The $H$-special radicals for $H$-module algebras}
J.R. Fisher \cite{Fi75} built up the general theory of
$H$-radicals for $H$-module
 algebras. We can easily give the definitions of the $H$-upper radical and
 the $H$-lower radical for $H$-module algebras as in \cite{Sz82}.
  In this section, we obtain
 some   properties of $H$-special radicals for $H$-module algebras.

\begin{Lemma}\label{9.1.1}
(1)   If $R$ is an $H$-module algebra and $E$ is a non-empty
subset of $R$, then
   $(E) = H\cdot E + R(H \cdot E) + (H\cdot E)R + R(H \cdot E)R$,
    where $(E)$ denotes the $H$-ideal generated by $E$ in  $R$.

(2) If $B$ is an $H$-ideal of $R$ and $C$ is an $H$-ideal of $B$,
then $(C)^3 \subseteq C,$  where $(C)$ denotes the $H$-ideal
generated by $C$ in $R$.
\end {Lemma}
{ \bf Proof. }  It is trivial. $\Box$

  \begin{Proposition}
  \label {9.1.2}   (1)  $R$ is $H$-semiprime iff $(H\cdot a)R(H\cdot a) = 0$ always implies
    $a = 0$ for any $a\in R$.

 (2) $R$ is $H$-prime iff $(H\cdot a)R(H\cdot b) = 0$ always implies
 $a = 0$  or $b = 0$ for any $a, b \in R$.
\end{Proposition}
{\bf Proof.}
 If $R$ is an $H$-prime module algebra and
$(H\cdot a)R(H\cdot b) = 0$ for $a, b \in R$, then $(a)^2 (b)^2 =
0$, where $(a)$ and $(b)$ are the $H$-ideals generated by  $a$ and
$b$  in $R$ respectively. Since $R$ is $H$-prime, $(a) = 0$ or
$(b)= 0$. Conversely, if $B$ and $C$ are  $H$-ideals of $R$ and $
BC = 0$, then $(H \cdot a)R(H\cdot b) = 0$  and $a = 0$ or $b = 0$
 for any $ a \in B, b \in C,$  which implies that $B = 0 $ or $C = 0$,
  i.e. $R$ is an $H$-prime module algebra.

Similarly, part (1) holds. $\Box$

 \begin{Proposition}\label {9.1.3}
  If $ I \lhd_{H} R$ and $I$ is an $H$-semiprime module algebra,
  then    \\
     (1)  $I \cap I^{*} = 0$;
     (2)  $I_{r} = I_{l} = I^{*}$;
     (3)  $I^{*} \lhd_H R$,
     where $I_{r} = \{ a \in R \mid  I(H\cdot a) = 0 \}$,
     $I_{l} = \{ a \in R \mid  (H\cdot a)I = 0 \}$,
     $I^{*} = \{ a \in R \mid  (H\cdot a)I = 0 = I(H\cdot a) \}.$
 \end{Proposition}
 {\bf Proof .} For any $x \in I^* \cap I$, we have that  $I(H\cdot x) = 0$
 and  $(H\cdot x)I(H\cdot x) = 0$. Since $I$ is an $H$-semiprime
 module algebra,  $x =0$,
 i.e. $I \cap I^* =0$.

 To show $I^* = I_r$,
 we only need to show that $(H \cdot x)I = 0$ for any $x \in I_r$.
 For any $y \in I, h \in H$, let $z = (h \cdot x)y$.
 It is clear that   $(H\cdot z)I(H\cdot z) = 0$.
 Since $I$ is an $H$-semiprime module algebra, $z =0$, i.e.
 $(H\cdot x)I = 0$. Thus $I^* = I_r$. Similarly, we can show that
 $I_l = I^*$.

 Obviously, $I^*$ is an ideal of $R$. For any $x \in I^*, h \in H$,
 we have $(H\cdot(h \cdot x)) I = 0 $. Thus $h\cdot x \in I^*$  by
 part (2), i.e. $I^*$ is an $H$-ideal of $R$.
 $\Box$

 \begin{Definition}\label {9.1.4} ${\cal K}$ is called an $H$-(weakly )special class if

     (S1)  ${\cal K}$ consists of $H$-(semiprime)prime module algebras.

     (S2)  For any $R \in {\cal K}$, if $0 \not= I \lhd_{H} R$  then
     $I \in {\cal K}$.

     (S3)  If $R$ is an $H$-module algebra and  $B\lhd_{H} R$ with
      $B \in {\cal K}$,
  then    $ R/B^{*}\in{\cal K}$,

where  $B^{*} = \{ a \in R\mid  (H\cdot a)B = 0 = B(H\cdot a)\}$.
    \end{Definition}

      It is clear that (S3) may be replaced by one of the following conditions:

     (S3')   If $B$ is an essential $H$-ideal of $R$(i.e. $B \cap I \not=0$
     for any non-zero $H$-ideal $I$ of $R$) and $B \in {\cal K}$, then
     $R\in {\cal K}$.

(S3")  If there exists an $H$-ideal $B$ of $R$ with $B^{*} = 0$
and $B\in{\cal K}$,
 then $R\in{\cal K}$.

 It is easy to check that if ${ \cal K}$ is an $H$-special class, then
 ${ \cal K}$ is an $H$-weakly special class.

 \begin{Theorem}\label {9.1.5}   If ${\cal K}$ is an $H$-weakly special class,
  then
   $r^{{\cal K}}(R) = \cap \{ I \lhd_{H} R \mid R/I \in {\cal K} \}$,
   where    $r^{{\cal K}}$ denotes the $H$-upper radical determined by $ {\cal K} $.
   \end{Theorem}
{\bf Proof.} If $I$ is a non-zero $H$-ideal of $R$ and $I \in
{\cal K}$,
 then $R/I^* \in {\cal K}$  by (S3) in Definition \ref{9.1.4}
 and  $I \not\subseteq I^* $ by Proposition  \ref {9.1.3}.
 Consequently, it follows from  \cite [Proposition 5]{Fi75} that
 $$r^{\cal K} (R) = \cap \{ I \mid I \hbox { \ is an \ } H \hbox {-ideal of } R
 \hbox { and } R/I \in {\cal K} \} \hbox { \ . \ \ }  \Box $$

    \begin{Definition}\label {9.1.6} If $r$ is a hereditary
    $H$-radical(i.e. if $R$ is
    an $r$-$H$-module algebra and   $B$ is an $H$-ideal of $R$, then so is $B$ )
    and any nilpotent $H$-module algebra is an $r$-$H$-module algebra,
   then $r$ is called a supernilpotent $H$-radical.
   \end{Definition}
\begin{Proposition}\label {9.1.7}  $r$ is a supernilpotent $H$-radical, then $r$
  is $H$-strongly hereditary,
  i.e. $r(I) = r(R) \cap I$   for any $I \lhd _{H} R$.
\end{Proposition}
{\bf Proof.} It follows from \cite [Proposition 4]{Fi75} . $\Box$

\begin{Theorem}\label {9.1.8}  If ${\cal K}$ is an $H$-weakly special class,
then     $r^{{\cal K}}$ is a supernilpotent $H$-radical.
\end{Theorem}
{\bf Proof.}     Let $r= r^{\cal K}.$
      Since every non-zero $H$-homomorphic image $R'$ of a nilpotent
$H$-module algebra $R$ is  nilpotent and is not $H$-semiprime, we
have that $R$ is an $r$-$H$-module algebra by Theorem  \ref
{9.1.5}. It remains to show that any $H$-ideal $I$ of
$r$-$H$-module
 algebra $R$ is an $r$-$H$-ideal.
If $I$ is not an $r$-$H$-module algebra,
 then there exists an $H$-ideal $J$ of $I$ such that
 $ 0 \not= I/J \in {\cal K}$. By (S3), $(R/J)/(I/J)^* \in {\cal K}$.
 Let
 $Q = \{x \in R \mid (H \cdot x)I
 \subseteq J$ and $I(H \cdot x) \subseteq J$ \}.
It is clear that
 $J$ and $Q$ are $H$-ideals of $R$ and  $Q/J= (I/J)^*$.
  Since $R/Q \cong (R/J)/(Q/J) = (R/J)/(I/J)^*$ and $R/Q$ is an $r$-$H$-
  module algebra, we have $(R/J)/(I/J)^*$ is an $r$-$H$-module algebra.
  Thus $R/Q = 0$ and $I^2 \subseteq J$, which  contradicts that
$I/J$ is a non-zero $H$-semiprime module algebra. Thus $I$ is an
$r$-$H$-ideal.
     $\Box$

\begin{Proposition}\label {9.1.9}  $R$ is $H$-semiprime iff for any  $0\not= a \in R,$
  there exists an $H$-$m$-sequence
  $\{ a_{n} \}$ in $R$ with $a_{7.1.1} = a$ such that $a_{n} \not= 0$ for all $n$.
  \end{Proposition}
{ \bf Proof.} If $R$ is $H$-semiprime, then for any $0 \not= a \in
R$, there exist
 $b_1 \in R$ , $h_1$ and $h_1' \in  H$  such that
 $ 0 \not= a_2 = (h_1 \cdot a_1) b_1(h_1' \cdot a_1)  \in
 (H \cdot a_1)R(H \cdot a_1)$ by Proposition \ref {9.1.2}, where $a_1 = a$.
 Similarly, for $0 \not= a_2 \in R$, there exist
 $b_2 \in R$ and $h_2$ and $h_2' \in  H$  such that
 $ 0 \not= a_3 = (h_2 \cdot a_2) b_2(h_2' \cdot a_2)  \in
 (H \cdot a_2)R(H \cdot a_2),$  which implies that
 there exists an $H$-$m$-sequence $\{ a_n \}$  such that
  $a_n \not= 0$ for any natural number $n$.  Conversely, it is trivial.
$\Box$

\section{$H$-Baer radical}

In this section, we give the characterization of $H$-Baer
radical($H$-prime radical) in many ways.

\begin{Theorem}\label {9.2.1} We define a property $r_{Hb}$ for
$H$-module algebras as follows: $R$ is an $r_{Hb}$-$H$-module
algebra iff every non-zero
   $H$-homomorphic image of $R$ contains a non-zero nilpotent $H$-ideal;
   then $r_{Hb}$ is an $H$-radical property.
   \end{Theorem}
  {\bf Proof.}   It is clear that every $H$-homomorphic image of
$r_{Hb}$-$H$-module algebra is an $r_{Hb}$-$H$-module algebra. If
every non-zero $H$-homomorphic image  $B$ of $H$-module algebra
$R$ contains a non-zero  $r_{Hb}$-$H$-ideal $I$, then $I$ contains
a non-zero nilpotent $H$-ideal $J$. It is clear that $(J)$ is a
non-zero nilpotent $H$-ideal of $B$, where (J) denotes the
$H$-ideal generated by $J$ in $B$. Thus $R$ is an
$r_{Hb}$-$H$-module algebra. Consequently, $r_{Hb}$ is an
$H$-radical property.
  $\Box$

$r_{Hb}$    is called $H$-prime radical or $H$-Baer radical.

  \begin{Theorem}\label {9.2.2}  Let
  $${\cal E} = \{ R \mid R \hbox { is a nilpotent } H\hbox {-module
   algebra } \},$$ then
   $r_{\cal E} = r_{Hb}$,  where    $r_{\cal E}$   denotes the $H$-lower radical determined
 by ${\cal E}$.
  \end{Theorem}
{\bf Proof.}
 If $R$  is an $r_{Hb}$-$H$-module algebra, then every non-zero
$H$-homomorphic image $B$ of $R$ contains a non-zero nilpotent
$H$-ideal $I$. By the definition of the lower $H$-radical, $I$ is
an $r_{\cal E}$-$H$-module algebra. Consequently, $R$ is an
$r_{{\cal E}}$-$H$-module algebra.
 Conversely, since every nilpotent $H$-module algebra is an
 $r_{Hb}$-$H$-module algebra, $r_{\cal E} \leq r_{Hb}$.
$\Box$

   \begin{Proposition}\label {9.2.3} R is $H$-semiprime if and only if $r_{Hb}(R) = 0$.
    \end{Proposition}
 {\bf Proof.}
 If $R$ is $H$-semiprime with $r_{Hb}(R) \not= 0$, then there exists
 a non-zero nilpotent $H$-ideal $I$ of $r_{Hb}(R)$.
 It is clear that $H$-ideal $(I)$,
 which the $H$-ideal  generated by $I$ in $R$,
 is  a non-zero nilpotent $H$-ideal of $R$. This contradicts that
 $R$ is  $H$-semiprime. Thus $r_{Hb}(R) =0$.
 Conversely, if $R$ is an $H$-module algebra with $r_{Hb}(R)= 0$ and
 there exists a non-zero nilpotent $H$-ideal $I$ of $R$,
 then $I \subseteq r_{Hb}(R)$.
 We get a contradiction. Thus $R$ is $H$-semiprime if
 $r_{Hb}(R)=0$.
 $\Box$

 \begin{Theorem}\label {9.2.4}  If  ${\cal K} = \{ R \mid\ R$ is an $H$-prime module algebra\},
 then $\cal K$ is an $H$-special class and
      $r_{Hb} = r^{{\cal K}}.$
       \end{Theorem}
{\bf Proof.}       Obviously, $(S1)$ holds.
 If $I$ is a non-zero $H$-ideal of an $H$-prime module
algebra $R$  and $BC = 0 $ for
 $H$-ideals $B$ and $C$  of $I$, then $(B)^3(C)^3 =0$
 where $(B)$ and $(C)$ denote the $H$-ideals generated by $B$ and $C$ in $R$
    respectively. Since $R$ is
 $H$-prime, $(B)= 0$ or $(C) =0$, i.e. $B =0$ or $C = 0$.
Consequently, $(S2)$ holds.
 Now we shows that (S3) holds. Let $B$ be an $H$-prime module algebra
 and be an $H$-ideal of $R$. If $ JI \subseteq B^*$ for $H$-ideals
 $I$ and $J$ of $R$, then $(BJ)(IB) = 0$, where
 $B^* = \{ x \in R \mid (H \cdot x)B = 0 = B(H \cdot x) \}$.
 Since $B$ is an $H$-prime module algebra, $BJ=0$ or $IB = 0$. Considering $I$ and $J$
 are $H$-ideals, we have that $B(H \cdot J)= 0$ or
 $(H \cdot I)B=0$. By Proposition \ref {9.1.3}, $J \subseteq B^*$
 or $I \subseteq B^*$, which implies that $R/B^* $  is an $H$-prime
 module algebra. Consequently, $(S3)$  holds and so ${\cal K}$ is an $H$-special
 class.

Next we show that $r_{Hb} = r^{\cal K}$. By Proposition \ref
{9.1.5}, $r^{\cal K}(R) = \cap \{ I \mid I$ is an $H$-ideal of $R$
and $R/I \in {\cal K}\}$. If $R$ is a nilpotent $H$-module
algebra, then $R$ is an $r^{\cal K}$-$H$-module algebra. It
follows from Theorem \ref {9.2.2} that $r_{Hb} \leq r^{\cal K}$.
Conversely, if $r_{Hb}(R)=0$, then $R$ is an $H$-semiprime module
algebra
 by Proposition \ref {9.2.3}.
 For any $0 \not= a \in R$, there exist
 $b_1 \in R$ , $h_1, h_1' \in  H$  such that
 $ 0 \not= a_2 = (h_1 \cdot a_1) b_1(h_1' \cdot a_1)  \in
 (H \cdot a_1)R(H \cdot a_1)$, where $a_1 = a$.
 Similarly, for $0 \not= a_2 \in R$, there exist
 $b_2 \in R$ and $h_2, h_2' \in  H$  such that
 $ 0 \not= a_3 = (h_2 \cdot a_2) b_2(h_2' \cdot a_2)  \in
 (H \cdot a_2)R(H \cdot a_2)$. Thus
 there exists an $H$-$m$-sequence $\{ a_n \}$  such that
  $a_n \not= 0$ for any natural number $n$. Let
  $${\cal F} = \{ I
  \mid I \hbox { is  an } H \hbox {-ideal  of } R \hbox { and }
  I \cap \{ a_1, a_2, \cdots \}
  = \emptyset \}.$$
 By Zorn's Lemma, there exists a
maximal element $P$  in ${\cal F}$. If $I$ and $J$ are $H$-ideals
of $R$ and $I \not\subseteq P$ and
 $J \not\subseteq P$, then there exist natural numbers
 $n$  and $m$ such that  $a_n \in I$  and
 $a_m \in J$. Since $0 \not= a_{n +m + 1}
 = (h_{n+m}\cdot a_{n+m})b_{n+m}(h'_{n+m}\cdot a_{n+m}) \in IJ$,
 which implies that
 $IJ \not\subseteq P$ and so $P$ is an $H$-prime ideal of $R$. Obviously, $a \not\in P,$
 which implies that $a \not\in r^{\cal K}(R)$
 and $r^{\cal K}(R)= 0$. Consequently,
  $r^{\cal K} = r_{Hb}$.  $\Box$

   \begin{Theorem}\label {9.2.5}  $r_{Hb}(R) = W_{H}(R)$.
     \end {Theorem}
{\bf Proof.} If $0 \not= a \not\in W_H(R)$, then there exists an
$H$-prime ideal $P$ such that $a \not\in P$  by the proof of
Thoerem \ref {9.2.4}. Thus $a \not\in r_{Hb}(R),$ which implies
that $r_{Hb}(R) \subseteq W_H(R)$.
          Conversely, for any $x \in W_H(R)$,
let $\bar R = R/r_{Hb}(R)$. Since $r_{Hb}(\bar R)=0$, $\bar R $ is
an
 $H$-semiprime module algebra by Proposition \ref {9.2.3}. By
 the proof of Theorem \ref {9.2.4}, $W_H(\bar R)=0$.
 For   an $H$-$m$-sequence $\{ \bar a_n \}$
 with $\bar a_1 = \bar x$ in  $\bar R$,  there exist
 $\overline b_n \in \overline R$ and $h_n, h_n' \in H$ such that
 $$ \overline a_{n+1} =
 (h_n \cdot \overline a_n) \overline b_n (h_n' \cdot \overline a_n)$$
                             for any natural number $n$.
Thus there exists
 $ a_n' \in  R$  such that  $a_1'=x $  and
 $  a_{n+1}' = (h_n \cdot  a_n')  b_n (h_n' \cdot  a_n')$ {~}{~}
                             for any natural number $n$.
 Since $ \{ a_n'\}$ is an $H$-$m$-sequence with $a_1' = x$ in $R$,
 there exists a natural number $k$ such that $a_k' =0$. It is easy
 to show that $\overline a_n= \overline a_n'$
  for any natural number $n$ by induction. Thus $ \bar a_k = 0$
  and $\overline x \in W_H(\overline R)$.
  Considering $W_H(\overline R)=0$, we have $x \in r_{Hb}(R)$ and
  $W_H(R) \subseteq r_{Hb}(R)$. Therefore
  $W_H(R) = r_{Hb}(R)$. $\Box$

 \begin{Definition}\label {9.2.6}  We define an $H$-ideal  $N_{\alpha}$ in $H$-module
 algebra  $R$ for every ordinal number $\alpha$ as follows:

 (i)  $N_{0} = 0$.

 Let us assume that $N_{\alpha}$ is already defined for $\alpha
\prec\beta$.

  (ii)  If $\beta = \alpha + 1, N_{\beta}/N_{\alpha}$ is the sum
   of all nilpotent  $H$-ideals  of $R/N_{\alpha}$

  (iii)  If $\beta$ is a limit ordinal number, $N_{\beta} =
  {\sum}_{\alpha\prec\beta} N_{\alpha}$.

  By set theory, there exists an ordinal number $\tau$ such that
  $N_{\tau} =  N_{\tau + 1}$.
   \end{Definition}

\begin{Theorem} \label {9.2.7} $N_{\tau} = r_{Hb}(R) =
\cap \{ I \mid I \hbox { is
   an } H \hbox {-semiprime ideal of } R \}$.
   \end{Theorem}
{\bf Proof.}
  Let $ D =\cap \{ I \mid I \hbox { is
   an } H \hbox {-semiprime ideal of } R \}.$
 Since  $R/N_{\tau}$ has not any non-zero nilpotent $H$-ideal, we have  that
$r_{Hb}(R) \subseteq N_{\tau}$ by Proposition \ref {9.2.3}.
Obviously, $D \subseteq r_{Hb}(R).$
 Using transfinite  induction, we can show that
 $N_{\alpha } \subseteq I$ for every $H$-semiprime ideal $I$ of $R$ and every
 ordinal number $\alpha $
 (see the proof of \cite [Theorem 3.7] {ZC91} ).
 Thus $ N_{\tau} \subseteq D,$  which completes the proof.
                        $\Box$

   \begin{Definition}
   \label {9.2.8} Let $ \emptyset \not= L \subseteq H$.
   An    $H$-$m$-sequence $\{ a_n \}$ in
    $R$ is called an $L$-$m$-sequence with beginning $a$ if
     $a_{7.1.1} = a$  and $a_{n+1} = (h_{n}.a_{n})b_{n}(h_{n}'.a_{n})$
    such that  $h_{n}, h_{n}'\in L$  for all $n$.
  For every $L$-$m$-sequence $\{a_{n} \}$ with $a_{7.1.1} = a$, there exists a
  natural number $k$ such that $a_{k} = 0$, then $R$ is called an $L$-$m$-nilpotent
  element,  written as  $W_{L}(R) = \{ a \in R \mid a$ is an $L$-$m$-nilpotent
  element\}.
     \end{Definition}

Similarly, we have
  \begin{Proposition}\label {9.2.9}  If $L \subseteq H$ and $H = kL$, then

(i)  $R$ is $H$-semiprime iff $(L.a)R(L.a) = 0$ always implies
  $ a = 0 $ for any  $a\in R$.

 (ii)    $R$ is $H$-prime iff $(L.a)R(L.b) = 0$ always implies
$ a = 0$  or $b = 0$       for any     $a, b \in R$.

(iii) $R$ is $H$-semiprime if and only if for any $0 \not= a \in
R$, there exists an $L$-$m$-sequence $\{a_{n} \}$  with $a_1 = a$
such that $a_{n}\not= 0 $ for all $n$.

(iv)  $W_{H}(R) = W_{L}(R)$.
\end{Proposition}

\section { The $H$-module theoretical characterization of $H$-special
radicals }
 If $V$ is an algebra over $k$ with unit and $x \otimes 1_V=0$  always implies
 that $x=0$ for any right $k$-module $M$ and for any $x \in M$,
 then $V$ is called a  faithful algebra to tensor. For example,
 if $k$ is a field,   then
 $V$ is  faithful to tensor for any algebra $V$ with unit.

In this section, we need to add the following condition:
  $H$ is  faithful  to tensor.

  We shall characterize
  $H$-Baer radical  $r_{Hb}$, $H$-locally nil radical $r_{Hl}$, $H$-Jacobson
 radical $r_{Hj}$ and $H$-Brown-McCoy radical $r_{Hbm}$ by $R$-$H$-modules.

 We can view  every $H$-module algebra $R$ as a sub-algebra of
 $R \# H$ since $H$ is   faithful to tensor.
 By computation, we have that $$ h\cdot a =
 \sum (1 \# h_1) a (1 \# S(h_2)) $$
 for any $h\in H, a \in R,$ where $S$ is the antipode of $H.$

  \begin{Definition}\label {9.3.1}
 An $R$-$H$-module $M$ is called an $R$-$H$-prime module
if for $M$ the following conditions are fulfilled:

(i)   $RM \not= 0$;

(ii)  If $x$ is an element of $M$ and $I$ is an $H$-ideal of $R$,
then $I(Hx) = 0$ always implies $x = 0$ or $I \subseteq
(0:M)_{R},$ where $(0:M)_R = \{ a \in R \mid  aM=0 \}$.
\end{Definition}
    \begin{Definition} \label {9.3.2}
   We associate to every $H$-module algebra $R$ a class ${\cal M}_{R}$ of
    $R$-$H$-modules. Then the class ${\cal M} = \cup {\cal M}_{R}$
is  called an $H$-special class of modules if the following
conditions are fulfilled:

(M1)  If $M \in {\cal M}_{R}$, then $M$ is an $R$-$H$-prime
module.

(M2)  If $I$ is an $H$-ideal of $R$ and $M \in {\cal M}_{I}$, then
$IM \in {\cal M}_{R}$.

(M3)  If $M \in {\cal M}_{R}$ and $I$ is an $H$-ideal of $R$ with
$IM \not= 0$, then $M \in {\cal M}_{I}$.

(M4)  Let $I$ be an $H$-ideal of $R$ and $\bar{R} = R/I$. If
 $M \in {\cal M}_{R}$
and $I \subseteq (0:M)_R$, then $M \in {\cal M}_{\bar {R}}$.
Conversely, if $M \in {\cal M}_{\bar {R}}$, then $M \in {\cal
M}_{R}$.
\end{Definition}

  Let ${\cal M}(R)$ denote $\cap \{ (0:M)_R \mid M \in {\cal M}_{R} \},$
  or $R$ when ${\cal M}_R = \emptyset $.
  \begin {Lemma} \label {9.3.3}
  (1)  If $M$ is an $R$-$H$-module, then $M$ is an $R\#H$-module.
In this case, $(0:M)_{R\#H}\cap R = (0:M)_{R}$ and $(0:M)_{R}$ is
an $H$-ideal of $R$;

 (2)  $R$ is a non-zero $H$-prime module algebra iff there exists a faithful
  $R$-$H$-prime
 module $M$;

(3) Let $I$ be an $H$-ideal of $R$ and $\bar R= R/I$. If $M$ is an
$R$-$H$-(resp. prime, irreducible)module and $I \subseteq (0:M)_R
$, then $M$ is an $\overline R$-$H$-(resp. prime,
irreducible)module (defined by $h \cdot (a + I) = h \cdot a$ and
$(a+I) x = ax )$. Conversely, if $M$ is an $\bar R$-$H$-(resp.
prime irreducible)module, then
 $M$ is an $R$-$H$-(resp. prime, irreducible)module(defined by $h \cdot a = h \cdot (a +I)$
 and $ax = (a+I)x$). In the both cases, it is always true that
$R/(0:M)_R \cong \overline R/(0:M)_{\overline R}$;

 (4)  $I$ is an $H$-prime ideal of $R$ with $I \not= R$ iff there exists an $R$-$H$-prime
 module $M$ such that $I=(0:M)_R$;

(5) If $I$ is an $H$-ideal of $R$ and $M$ is an $I$-$H$-prime
module, then $IM$ is an $R$-$H$-prime module with $ (0:M)_I =
(0:IM)_R\cap I$;

(6) If $M$ is an $R$-$H$-prime module and $I$ is an $H$-ideal of
$R$ with $IM \not=0$, then $M$ is an $I$-$H$-prime module;

(7) If $R$ is an $H$-semiprime module algebra with  one side unit,
then $R$ has a unit.
 \end{Lemma}
 {\bf Proof}. (1)  Obviously, $(0:M)_R = (0:M)_{R\#H} \cap R$. For any $h \in H,
  a \in (0:M)_R$,
 we see that $(h\cdot a)M = \sum (1\# h_1) a ( 1 \#S(h_2))M
 \subseteq \sum (1 \# h_1)aM=0$ for any $h\in H, a \in R$. Thus $h\cdot a \in (0:M)_R$,
 which implies $(0:M)_R$ is an $H$-ideal of $R$.

 (2) If $R$ is an $H$-prime module algebra, view  $M = R$ as an
$R$-$H$-module. Obviously, $M$ is faithful. If $I(H\cdot x) = 0$
 for $0 \not= x \in M $
 and an $H$-ideal $I$ of $R$, then
$ I(x)= 0$  and $I = 0$, where (x) denotes the $H$-ideal generated
by $x$  in $R$. Consequently,  $M$ is a faithful $R$-$H$-prime
module. Conversely, let $M$ be a faithful $R$-$H$-prime module. If
$IJ = 0$ for two $H$-ideals  $I$ and $J$ of $R$ with $J \not= 0$,
then $JM \not=0$ and there exists $0 \not= x \in JM$ such that
$I(H x) = 0$. Since $M$ is a faithful  $R$-$H$-prime module, $I =
0$. Consequently,  $R$ is $H$-prime.

(3)  If $M$ is an $R$-$H$-module, then it is clear that $M$ is a
(left)$\overline R$-module and $h  (\overline ax) = h  ( a x) =
  \sum (h_1 \cdot  a )(h_2 x)
  = \sum \overline{ (h_1 \cdot a) }(h_2  x)
  = \sum (h_1 \cdot \overline a)(h_2  x)$ {~}{~}
  for any $h \in H$,  $a \in R$ and $x \in M$. Thus
  $M$ is an $\overline R$-$H$-module. Conversely, if $M$ is an $\overline R$-$H$-module,
  then $M$ is an (left) $R$-module and
 $$h  (ax)
= h (\overline a x)  =
  \sum (h_1 \cdot \overline a )(h_2 x)
  = \sum \overline{h_1 \cdot a}(h_2  x)
  = \sum (h_1 \cdot  a)(h_2  x)$$
    {~}{~} for any $h \in H$, $a \in R$
  and $x \in M$.  This shows that
  $M$ is an $R$-$H$-module.

  Let  $M$ be an $R$-$H$-prime module and $I$ be an $H$-ideal of $R$
  with $I \subseteq (0:M)_R$. If $\overline J(Hx) =0$ for $0 \not= x \in M$
   and an $H$-ideal $J$ of $R$, then $J(Hx)=0$  and $J \subseteq (0:M)_R$.
   This shows that $\bar J \subseteq (0:M)_{ \overline R}$.
   Thus $M$ is an $R$-$H$-prime module. Similarly,
   we can show the other assert.

(4)  If $I$ is an $H$-prime ideal of $R$ with $R\not=I$,
 then $\overline R = R/I$ is an $H$-prime
module algebra. By Part (2), there exists a faithful $\overline
R$-$H$-prime module $M$. By part (3), $M$ is an $R$-$H$-prime
module with $(0:M)_R = I$. Conversely, if there exists a
$R$-$H$-prime $M$ with $I = (0:M)_R$, then $M$ is a faithful
$\overline R$-$H$-prime module by part (3) and $I$ is an $H$-prime
ideal of $R$ by part (2).

(5)  First, we show that $IM$ is an $R$-module. We define
 \begin {eqnarray}  \label {e2.3.3.1} a (\sum_{i}a_i x_i )= \sum_{i} (a a_i)x_i
\end {eqnarray}
for any $a \in R$ and $\sum_i a_ix_i \in IM$, where $a_i \in I$
and $x_i \in M$. If $\sum_{i}a_i x_i = \sum_{i}  a_i' x_i'$ with
$a_i$, $a_i'  \in  R$, $x_i$,  $x_i' \in M$, let $y =
\sum_{i}(aa_i )x_i - \sum_{i}  (aa_i') x_i'$. For any $b \in I$
and $h \in H$, we see that
\begin {eqnarray*}
b(h y) &=& \sum_{i} b\{h  [(aa_i) x_i -  (aa_i') x_i']\} \\
 &=& \sum_{i}\sum_{(h)} b\{[(h_1 \cdot (a a_i)]( h_2  x_i) -
 [h_1 \cdot (a a_i')] (h_2 x_i') \}  \\
  &=& \sum_{(h)} \sum_i \{b[(h_1 \cdot a)(h_2 \cdot a_i)](h_3 x_i)) -
 b[(h_1 \cdot a)(h_2 \cdot a_i')] (h_3  x_i')\}  \\
  &=& \sum_{(h)} \sum_{i} b(h_1 \cdot a)[h_2  (a_i x_i) -
 h_2 (a_i'x_i')] \\
 &=& \sum_{(h)}  b(h_1 \cdot a)h_2  \sum_i [a_i x_i -
a_i'x_i'] = 0.
 \end {eqnarray*}
 Thus  $I(H y)=0$.
 Since $M$ is an $I$-$H$-prime module and $IM \not= 0$,
 we have that $y=0$. Thus
 this definition in (\ref {e2.3.3.1}) is well-defined.
 It is easy to check that $IM$ is an $R$-module.
 We see that
\begin {eqnarray*}
h (a \sum_i a_ix_i) &=& \sum_i h[(aa_i)x_i]   \\
 &=& \sum_i \sum_h [h_1 \cdot (aa_i)][h_2 x_i]  \\
 &=&  \sum_i \sum_h [(h_1 \cdot a)(h_2 \cdot a_i)](h_3 x) \\
 &=& \sum_h (h_1 \cdot a) \sum_i (h_2 \cdot a_i)(h_3 x_i) \\
  &=& \sum_h (h_1 \cdot a) [ h_2 \sum_i (a_ix_i)]
\end {eqnarray*}
for any $h \in H$ and $\sum_i a_ix_i \in IM.$ Thus  $IM$ is an
$R$-$H$-module.

 Next, we show that $(0:M)_I = (0:IM)_R \cap I$.
 If $a \in (0:M)_I$, then $aM =0$  and $aIM=0$,
 i.e. $a \in (0:IM)_A \cap I$. Conversely, if $a \in (0:IM)_R \cap I$,
 then $aIM=0$. By part (1), $(0:IM)_R$ is an $H$-ideal of $R$.
 Thus $(H\cdot a)IM =0$ and $(H \cdot a)I \subseteq (0:M)_I$.
 Since $(0:M)_I$ is an $H$-prime ideal of $I$ by part (4),
 $a \in (0:M)_I$. Consequently,  $(0:M)_I = (0:IM)_R \cap I$.

Finally, we show that $IM$ is an $R$-$H$-prime module. If $RIM=0$,
then $RI \subseteq (0:M)_R$ and  $I \subseteq (0:M)_R,$ which
contradicts that $M$ is an $I$-$H$-prime module. Thus $RIM
\not=0$. If $J(Hx) = 0$ for $ 0\not= x \in IM $  and an $H$-ideal
$J$  of $R$, then $JI(Hx) \subseteq J(Hx) = 0$. Since $M$ is an
$I$-$H$-prime module, $JI \subseteq (0:M)_I$ and $J(IM) =0$.
Consequently,  $IM$ is an $R$-$H$-prime module.

(6)  Obviously, $M$ is an $I$-$H$-module. If $J(Hx) =0$ for $0
\not= x \in M$ and an $H$-ideal $J$ of $I$, then $(J)^3(H x) = 0$
and $(J)^3 \subseteq (0:M)_R$, where (J) denotes the $H$-ideal
generated by $J$ in $R$.
 Since $(0:M)_R$  is an $H$-prime ideal of $R$, $(J) \subseteq (0:M)_R$
  and $J \subseteq (0:M)_I$. Consequently,  $M$ is an $I$-$H$-prime module.

   (7) We can assume that $u$ is a right unit of $R$.
   We see that
   $$(h \cdot (au-a))b = \sum (1 \# h_1 ) (au-a)(1 \# S(h_2)) b =0$$
   for any $a, b \in R, h \in H.$  Therefore
   $(H\cdot (au-a)) R=0$ and $au =a$, which implies that $R$ has a unit.
     $\Box$

 \begin{Theorem} \label {9.3.4}
 (1)  If ${\cal M}$ is an $H$-special class of modules and ${\cal K}$
= \{ $R \mid$ there exists a faithful $R$-$H$-module $M \in {\cal
M}_{R}$\}, then ${\cal K}$ is an $H$-special class and $r^{{\cal
K}}(R) = {\cal M}(R)$.

(2)  If ${\cal K}$ is an $H$-special class and ${\cal M}_{R}$ = \{
$M \mid M$ is an $R$-$H$-prime module and $R/(0:M)_R \in {\cal
K}$\}, then ${\cal M} = \cup {\cal M}_{R}$ is an $H$-special class
of modules and $r^{{\cal K}}(R) = {{\cal M}}(R)$.
 \end{Theorem}
 { \bf Proof.}
 (1) By Lemma \ref {9.3.3}(2), $(S1)$ is satisfied. If $I$ is a non-zero
$H$-ideal of $R$  and $R \in {\cal K}$, then there exists a
faithful $R$-$H$-prime module  $M \in {\cal M}_R$. Since $M$ is
faithful, $IM \not=0$ and $M \in {\cal M}_I$ with $(0:M)_I =
(0:M)_R\cap I=0$ by $(M3)$. Thus $I \in {\cal K}$ and $(S2)$ is
satisfied. Now we show that $(S3)$ holds. If $I$ is an $H$-ideal
of $R$ with $I \in {\cal K}$, then there exists a faithful
$I$-$H$-prime module $M\in {\cal M}_I$. By $(M2)$ and Lemma \ref
{9.3.3}(5),
 $IM\in {\cal M}_R$ and $0=(0:M)_I=(0:IM)_R \cap I$.
Thus $(0:IM)_R \subseteq I^* $. Obviously, $I^* \subseteq
(0:IM)_R$. Thus $I^* = (0:IM)_R$.
 Using $(M4)$, we have that $IM \in {\cal M}_{\overline R}$
 and $IM$ is a faithful $\overline R$-$H$-module with $\overline R = R/I^*$.
 Thus $R/I^* \in {\cal K}$. Therefore ${\cal K}$ is an $H$-special class.

It is clear that
\begin {eqnarray*}
\{ I \mid I \hbox { is an } H \hbox {-ideal of } R \hbox { and }
R/I \in {\cal K} \}  &=& \{(0:M)_R \mid M \in {\cal M}_R  \}.
  \end {eqnarray*}
  Thus $r^{\cal K}(R) = {\cal M}(R)$.

 (2) It is clear that $(M1)$ is satisfied. If $I$ is an $H$-ideal of
$R$ with $M \in {\cal M}_I$, then $M$ is an $I$-$H$-prime module
with $I/(0:M)_I\in {\cal K}$. By Lemma \ref {9.3.3}(5), $IM$ is an
$R$-$H$-prime module with $(0:M)_I=(0:IM)_R \cap I$. It is clear
that
$$(0:IM)_R = \{ a \in R \mid (H \cdot a )I \subseteq (0:M)_I
\hbox { and } I(H \cdot a)  \subseteq (0:M)_I \} $$ and
$$(0:IM)_R/(0:M)_I= (I/(0:M)_I)^*.$$
Thus $R/(0:IM)_R \cong (R/(0:M)_I)/((0:IM)_R/(0:M)_I) =
(R/(0:M)_I)/(I/(0:M)_I)^* \in {\cal K},$ which implies that $IM
\in {\cal M}_R$ and $(M2)$ holds.
 Let $M \in {\cal M}_R$ and $I$ be an $H$-ideal of $R$
 with $IM \not=0$. By Lemma \ref {9.3.3}(6), $M$ is an $I$-$H$-prime
 module and
 $I/(0:M)_I = I/((0:M)_R \cap I) \cong (I + (0:M)_R)/(0:M)_R$.
 Since $R/(0:M)_R \in {\cal K}$, $I/(0:M)_I \in {\cal K}$ and
 $M \in {\cal M}_I$.
 Thus $(M3)$ holds.
 It follows from Lemma \ref{9.3.3}(3) that $(M4)$ holds.

It is clear that
 $$\{I \mid I \hbox { is an } H \hbox{-ideal of } R \hbox { and }
 0 \not= R/I \in {\cal K} \}
 = \{ (0:M)_R \mid M \in {\cal M}_R \}.$$
 Thus  $r^{\cal K}(R) = {\cal M}(R)$.
$\Box $

  \begin{Theorem} \label {9.3.5}
  Let  ${\cal M}_{R}$ =\{ $M \mid M$
  is an  $R$-$H$-prime module\} for any $H$-module algebra $R$
   and ${\cal M} = \cup {\cal M}_{R}$.
   Then ${\cal M}$ is an $H$-special class of modules and
  ${\cal M}(R) = r_{Hb}(R)$.
  \end{Theorem}
{\bf Proof.} It follows from Lemma \ref {9.3.3}(3)(5)(6) that
  ${\cal M}$ is an $H$-special class of modules. By Lemma \ref {9.3.3}(2),
    $$\{ R \mid R  \hbox { is an } H \hbox {-prime module algebra with }
    R \not=0 \} = $$
    $$\{   R \mid \hbox { there exists a faithful } R \hbox {-}H
   \hbox {-prime module }\}. $$ \\
   Thus  $r_{Hb}(R)= {\cal M}(R)$ by Theorem \ref {9.2.4}(1).
$\Box$

  \begin{Theorem} \label {9.3.6}
  Let  ${\cal M}_{R}$ =\{ $M \mid M$
  is an  $R$-$H$-irreducible module\} for any $H$-module algebra $R$
   and ${\cal M} = \cup {\cal M}_{R}.$
  Then ${\cal M}$ is an $H$-special class of modules and
  ${\cal M}(R) = r_{Hj}(R),$            where $r_{Hj}$ is the $H$-Jacobson
   radical of $R$ defined in \cite {Fi75}.
    \end{Theorem}

{\bf Proof.}                       If $M$ is an
$R$-$H$-irreducible module and $J (Hx) = 0$  for $0 \not= x \in M$
and an $H$-ideal $J$ of $R$, let  $N= \{ m \in M \mid J(H  m) = 0
\}$. Since $J(h (am)) = J(\sum_h (h_1 \cdot a)(h_2  m)) = 0$, $am
\in N$ for any $m \in N, h \in H, a \in R,$  we have
  that $N$ is an $R$-submodule of $M$.
Obviously, $N$ is an $H$-submodule of $M$. Thus $N$ is an
$R$-$H$-submodule of $M$. Since $N \not= 0$, we have that $N = M$
and $JM= 0$, i.e. $J \subseteq (0:M)_R$. Thus $M$ is an
$R$-$H$-prime  module and (M1) is satisfied. If $M$ is an
$I$-$H$-irreducible module and $I$ is an $H$-ideal, then $IM$ is
an $R$-$H$-module. If $N$ is an $R$-$H$-submodule of $IM$, then
$N$ is also an $I$-$H$-submodule of $M,$ which implies that $N=0$
or $N=M$. Thus $(M2)$ is satisfied. If $M$ is an
$R$-$H$-irreducible module and $I$ is an $H$-ideal of $R$ with
$IM\not=0$, then $IM=M$. If $N$ is an non-zero $I$-$H$-submodule
of $M$, then $IN$ is an $R$-$H$-submodule of $M$ by Lemma
\ref{9.3.3}(5) and $IN = 0 $  or $IN=M$. If $IN=0$, then $I
\subseteq (0:M)_R$ by the above proof and $IM=0$. We get a
contradiction. If $IN=M$, then $N=M$. Thus $M$ is an
$I$-$H$-irreducible module and $(M3)$ is satisfied.

It follows from Lemma \ref {9.3.3}(3) that (M4) holds. By Theorem
\ref {9.3.4}(1), ${\cal M}(R)= r_{Hj}(R)$. $\Box$

 J.R. Fisher
  \cite[Proposition 2]{Fi75} constructed
  an $H$-radical $r_{H}$ by a common hereditary radical $r$ for algebras,
  i.e. $r_H(R) = (r(R):H) =
  \{ a \in R \mid h\cdot a \in r(R)  \hbox { for any }
  h \in H \}.$  Thus we can get
  $H$-radicals $r_{bH}, r_{lH}, r_{jH}, r_{bmH}$.

  \begin{Definition} \label {9.3.7}
  An $R$-$H$-module $M$ is called an $R$-$H$-$BM$-module, if for
    $M$ the following conditions are fulfilled:

    (i)  $RM \not=0$;

    (ii)  If $I$ is an $H$-ideal of $R$ and $I \not\subseteq (0:M)_R$, then there exists
    an element $u \in I$ such that $m = um$ for all $m  \in M$.
\end{Definition}

\begin{Theorem} \label {9.3.8}
  Let ${\cal M}_{R}$ = \{ $M \mid M$ is an $R$-$H$-$BM$-module\} for
  every $H$-module algebra $R$ and ${\cal M} = \cup {\cal M}_{R}$. Then ${\cal M}$ is an
  $H$-special class of  modules.
\end{Theorem}
{\bf Proof.} It is clear that {\cal M} satisfies $(M_1)$  and
$(M_4)$. To prove $(M_2)$ we exhibit: if $I \lhd _H R$  and $M \in
{\cal M}_I$, then $M$ is an $I$-$H$-prime module and $IM$ is an
$R$-$H$-prime module. If $J $ is an $H$-ideal of $R$ with $J
\not\subseteq (0:M)_R, $  then $JI$ is an $H$-ideal of $I$ with
$JI \not\subseteq (0:M)_I.$  Thus there exists an element $u \in
JI \subseteq J$ such that
 $um = m $ for every $m \in M$. Hence
$IM \in {\cal M}_R.$

To prove $(M_3)$, we exhibit: if $M \in {\cal M}_R$ and $I$ is an
$H$-ideal of $R$ with $IM \not=0.$ If $J$ is an $H$-ideal of $I$
with $J \not\subseteq (0:M)_I,$ then $ (J) \not\subseteq (0:M)_R,$
where $ (J)$  is the $H$-ideal generated by $J$ in $R.$  Thus
there exists  an elements $u\in  (J)$ such that  $um=m$  for every
$m \in M.$ Moreover, $$m=um=uum=uuum = u^3 m $$ and $u^3 \in J.$
Thus $M \in {\cal M}_I.$ $\Box$

\begin{Proposition} \label {9.3.9}
  If $M$ is an $R$-$H$-$BM$-module, then
  $R/(0:M)_R$ is an  $H$-simple module algebra with unit.
\end {Proposition}

{\bf Proof. }  Let $I$ be any $H$- ideal of $R$  with $I
\not\subseteq (0:M)_R.$ Since $M$ is an $R$-$H$-$BM$-module, there
exists an element $u \in I$ such that $uam =am$  for every $m\in
M, a \in R.$ It follows that $a -ua \in (0:M)_R,$ whence
$R=I+(0:M)_R.$  Thus $(0:M)_R$ is a maximal $H$-ideal of $R.$
Therefore $R/(0:M)_R$  is an $H$-simple module algebra.

Next we shall show that $R/(0:M)_R$ has a unit. Now $R
\not\subseteq (0:M)_R,$ since $RM \not=0.$  By the above proof,
there exists an element $u \in R$ such that $a-ua \in (0:M)_R$ for
any $a \in R.$  Hence $R/(0:M)_R$ has a left unit. Furthermore,
by Lemma \ref {9.3.7} (7) it has a unity element. $\Box$

\begin{Proposition} \label {9.3.10}
  If $R$ is an $H$-simple-module algebra with unit, then
  there exists a faithful $R$-$H$-$BM$-module.
  \end {Proposition}
{\bf Proof. }  Let $M= R.$ It is clear that $M$ is a faithful
$R$-$H$-$BM$- module. $\Box$

\begin{Theorem} \label {9.3.11}
  Let ${\cal M}_{R}$ = \{ $M \mid M$ is an $R$-$H$-$BM$-module\} for
  every $H$-module algebra $R$ and ${\cal M} = \cup {\cal M}_{R}$. Then
$r_{Hbm}(R) = {\cal M}(R)$,
  where $r_{Hbm}$ denotes the $H$-upper radical determined by $\{R \mid R$ is
  an $H$-simple     module algebra with unit \}.
\end{Theorem}
 {\bf Proof.} By Theorem \ref {9.3.8}, ${\cal M}$ is an $H$-special class
 of modules. Let $${\cal K} = \{   R \mid \hbox { there exists a faithful }
 R  \hbox {-}H \hbox {-}BM\hbox {-module }   \}.$$  By Theorem \ref {9.3.4}(1),
   ${\cal K}$ is an $H$-special class  and  $r ^{\cal K}(R) = {\cal M}(R).$
  Using Proposition \ref  {9.3.9} and \ref {9.3.10}, we have that

$$ {\cal K} = \{  R \mid R \hbox { is an $H$-simple module algebra with
unit }\}.$$ Therefore ${\cal M}(R) = r_{Hbm}(R).$  $\Box$

  Assume that
$H$ is a finite-dimensional semisimple Hopf
  algebra with
 $t \in \int_{H}^{l} $ and $\epsilon(t) = 1$. Let
 $$G_{t}(a) =
  \{z \mid z = x + (t.a)x + \sum(x_{i}(t.a)y_{i} + x_{i}y_{i})
  \hbox { \ for all \ } x_{i}, y_{i}, x \in R \}.$$
   $R$ is called an $r_{gt}$-$H$-module algebra, if $a \in G_{t}(a)$ for all
    $a \in R$.

\begin{Theorem}\label {9.3.12}
$r_{gr}$ is an $H$-radical property of $H$-module algebra and
$r_{gt} =r_{Hbm}$.
\end{Theorem}

{\bf Proof.}  It is clear that any $H$-homomorphic image of
$r_{gt}$-$H$- module algebra is an $r_{gt}$-$H$-module algebra.
Let $$N= \sum \{I \lhd _H \mid I \hbox { is an }  r_ {gt} \hbox
{-} H \hbox {-ideal of } R\}.$$ Now we show that $N$ is an
$r_{gt}$-$H$-ideal of $R$. In fact, we only need to show  that
$I_1 +I_2$ is an $r_{gt}$-$H$-ideal  for any two
$r_{gt}$-$H$-ideals  $I_1$ and $I_2$. For any $a \in I_1, b \in
I_2$, there exist $x,  x_i, y_i \in R$ such that
$$a = x + (t \cdot a)x + \sum_i (x_i (t \cdot a)y_i + x_iy_i).$$
Let $$c = x + (t \cdot (a +b))x +   \sum x_i (t \cdot (a +b))y_i +
x_iy_i  \hbox { \ \ \ } \in G_t(a +b).$$ Obviously, $$ a +b -c = b
- (t \cdot b)x - \sum x_i (t\cdot b)y _i \hbox { \ \  \ } \in
I_2.$$ Thus there exist $w, u_j, v_j \in R$ such that
$$a+b-c= w + (t \cdot (a+b-c))w + \sum _j (u_j(t \cdot (a+b-c))v_j + u_jv_j).$$
Let $d= (t\cdot (a +b))w +w + \sum _j (u_j(t \cdot (a +b))v_j
+u_jv_j)$ and $e = c - \sum _j u_j(t \cdot c)v_j - (t \cdot c)w.$
By computation, we have that $$a +b =d +e.$$ Since $c \in G_t (a
+b)$ and $d \in G_t(a+b)$, we get that $e \in G_(a +b)$ and $a
+b\in G_t(a+b),$ which implies that $I_1 +I_2$
 is an $r_{gt}$-$H$-ideal.

Let $\bar R = R/N$ and $\bar B $ be an $r_{gt}$-$H$-ideal of $\bar
R.$ For any $a \in B$, there exist $ x , x _i , y_i \in R$ such
that $$\bar a = \bar x+ (t \cdot \bar a) \bar x + \sum (\bar x_i
(t \cdot \bar a ) \bar y_i + \bar x_i \bar  y_i)$$ and  $$x + (t
\cdot  a)  x + \sum ( x_i (t \cdot  a )  y_i +  x_i   y_i) -a \in
N. $$ Let $$c = x + (t \cdot a)x + \sum (x_i (t \cdot a) y _i +
x_i y_i) \in G_t(a).$$ Thus there exist $w, u_j, v_j \in R$ such
that
$$a -c = (t\cdot (a-c))w +w + \sum (u_j (t \cdot (a-c))v_j + u_jv_j) $$ and
$$a = (t\cdot a)w +w + \sum u_j (t \cdot a)v_j + u_jv_j +c
-(t \cdot c)w - \sum u_j(t \cdot c)v_j \hbox { \ \  \ }\in G_t(a)
,$$ which implies that $B$ is an $r_{gt}$-$H$-ideal and $\bar
B=0.$ Therefore $r_{gt}$ is an $H$ -radical property.  $\Box$

\begin{Proposition}\label {9.3.13}
If $R$ is an $H$-simple module algebra, then $r_{gr}(R) =0$ iff
$R$ has a unit.
\end{Proposition}
{\bf Proof.} If $R$ is an $H$-simple module algebra with unit $1$,
then $-1 \not\in G_t(-1)$ since
$$x + (t \cdot (-1))x + \sum (x_i (t \cdot (-1))y_i + x_iy_i )=0 $$
for any $x, x_i, y_i \in R$. Thus $R$ is $r_{gt}$-$H$-semisimple.
Conversely, if $r_{gt}(R)=0,$ then there exists $0 \not=a \not\in
G_t(a)$ and $G_t(a)=0,$ which implies that $ax +x=0$  for any $x
\in R.$  It follows from Lemma {9.3.3} (7) that $R$ has a unit.
$\Box$

\begin{Theorem}\label {9.3.14}
 $r_{gt} =r_{Hbm}$.
 \end{Theorem}
{\bf Proof.}  By Proposition \ref {9.3.13}, $r_{gt}(R) \subseteq
r_{Hbm}(R) $ for any $H$-module algebra $R.$  It  remains
 to show that if $a \not\in r_{gt}(R)$ then
$a \not\in r_{Hbm}(R) $. Obviously, there exists $b \in (a)$ such
that $b \not\in G_t(b),$  where $(a)$  denotes the $H$-ideal
generated by $a$ in $R.$
 Let $${\cal E} = \{ I \lhd _H R \mid G_t(b) \subseteq I , b \not\in I \}.$$
 By Zorn's Lemma, there exists a maximal element $P$ in ${ \cal E }$. $P$ is a
 maximal $H$-ideal of $R$, for, if $Q$  is an $H$-ideal of $R$  with
 $P \subseteq Q$ and $P\not= Q,$  then $b \in Q$   and
 $x = -bx + (bx +x) \in Q$  for any $x \in R.$ Consequently,
 $R/P$ is an $H$-simple module algebra with $r_{gt}(R/P)=0$.
 It follows from Proposition \ref {9.3.13}  that
 $R/P$  is an $H$-simple module algebra with unit and $r_{Hbm}(R) \subseteq P.$
 Therefore  $b \not\in r_{Hbm}(R)$  and so $a \not\in r_{Hbm}(R).$
 $\Box$

                                 \begin{Definition} \label {9.3.15}
    Let $I$ be an $H$-ideal of $H$-module algebra $R$, $N$ be an $R$-$H$-submodule
    of $R$-$H$-module $M$. $(N,I)$ are said to have ``L-condition'',
    if for any finite subset $F \subseteq I$, there exists a positive integer $k$
    such that $F^{k}N = 0$.
 \end{Definition}
  \begin{Definition} \label {9.3.16}
  An $R$-$H$-module $M$ is called an $R$-$H$-$L$-module, if for
   $ M$ the following conditions are fulfilled:

    (i)  $RM \not= 0$.

    (ii)  For every non-zero $R$-$H$-submodule $N$ of $M$ and every $H$-ideal $I$ of
    $R$, if $(N,I)$ has ``$L$-condition'', then $I \subseteq (0:M)_R$.
 \end{Definition}

    \begin{Proposition} \label {9.3.17}
If $M$ is an $R$-$H$-$L$-module, then $R/(0:M)_R$ is an
$r_{lH}$-$H$-semisimple and $H$-prime module algebra.
     \end {Proposition}
{\bf Proof.}  If $M$ is an $R$-$H$-$L$-module, let $\bar R =
R/(0:M)_R.$ Obviously, $\bar R $ is $H$-prime.
 If $\bar B$ is an   $r_{lH}$-$H$-ideal of $\bar R$, then
 $(M, B)$  has "$L$-condition" in $R$-$H$-module $M$, since for any
 finite subset $F$ of $B$, there exists a natural number $n$ such that
 $F^n \subseteq (0:M)_R$ and $F^nM=0.$  Consequently,
 $B \subseteq (0:M)_R$  and $\bar R$ is  $r_{lH}$-semisimple.
 $\Box$

    \begin{Proposition} \label {9.3.18}
    $R$ is a non-zero $r_{lH}$-$H$-semisimple and $H$-prime module algebra iff
    there exists a faithful $R$-$H$-$L$-module.
    \end {Proposition}
 {\bf Proof.}  If    $R$ is a non-zero $r_{lH}$-$H$-semisimple and
 $H$-prime module algebra, let $M = R$. Since $R$ is an $H$-prime module
 algebra, $(0:M)_R=0.$
 If $(N,B)$ has "$L$-condition" for non-zero $R$-$H$-submodule of $M$ and
 $H$-ideal $B$, then, for any finite subset $F$ of $B$, there exists an natural number
 $n$, such that  $F^n N =0$ and $F^n (NR)=0$, which implies that
 $F^n =0$ and $B$ is an $r_{lH}$-$H$-ideal, i.e.
 $B=0 \subseteq (0:M)_R$. Consequently,  $M$ is a faithful $R$-$H$-$L$-
 module.

 Conversely, if $M$ is a faithful $R$-$H$-$L$-module, then $R$ is an $H$-prime
 module algebra. If $I$ is an $r_{lH}$-$H$-ideal of $R$, then
 $(M,I)$ has ``$L$-condition", which implies $I=0 $  and
 $R$ is an $r_{lH}$-$H$-semisimple module algebra.
$\Box$

 \begin{Theorem} \label {9.3.19}
 Let ${\cal M}_{R}$ =
    \{ $M \mid M$ is an $R$-$H$-$L$-module\} for any $H$-module algebra
    $R$ and ${\cal M} = \cup {\cal M}_{R}$. Then
    ${\cal M}$ is an $H$-special class of modules
   and      ${\cal M} (R) = r_{Hl}(R),$
   where ${\cal K} = \{ R \mid R \hbox { is an } H \hbox{-prime module
 algebra with }  r_{lH}(R)=0 \}$  and
  $r_{Hl}= r^{\cal K}$.
  \end {Theorem}

{\bf Proof.}          Obviously, $(M1)$  holds. To show  that
$(M2)$ holds, we only need to show that if $I$ is an $H$-ideal of
$R$ and $M \in {\cal M}_I,$ then $IM \in {\cal M}_R.$ By Lemma
\ref {9.3.3}(5), $IM$ is an $R$-$H$-prime module. If $(N,B)$  has
the "$L$-condition"  for non-zero  $R$-$H$-submodule $N$ of $IM$
and $H$-ideal $B$ of $R$, i.e. for any finite subset $F$  of $B$,
there exists a natural number $n$ such that $F^n N=0$, then $(N,
BI)$ has "$L$-condition" in $I$-$H$-module $M$. Thus $BI \subseteq
(0:M)_I = (0:IM)_R \cap I$. Considering  $(0:IM)_R$  is an
$H$-prime ideal of $R$, we have that $B \subseteq (0:IM)_R$ or $I
\subseteq (0:IM)_R.$  If $I \subseteq (0:IM)_R$, then $I^2
\subseteq (0:M)_I$ and $I \subseteq (0:M)_I,$ which contradicts
$IM \not=0.$ Therefore $B \subseteq (0:IM)_R$ and so $IM$   is an
$R$-$H$-$L$- module.

To show that $(M3)$ holds, we only need to show that if $M \in
{\cal M}_R$ and $I \lhd _H R$ with $IM \not=0$, then $M\in {\cal
M}_I.$  By Lemma \ref {9.3.3}(6), $M$ is an $I$-$H$-prime module.
If $(N,B)$  has  the "$L$-condition"  for non-zero
$I$-$H$-submodule $N$ of $M$ and $H$-ideal $B$ of $I,$ then
 $IN$ is an $R$-$H$-prime module and $(IN, (B))$ has "$L$-condition"
in $R$-$H$-module $M,$ since for any finite subset $F$  of $(B)$,
$F^3 \subseteq B$ and  there exists a natural number $n$ such that
$F^{3n}I N\subseteq F^{3n}N=0$, where $(B)$ is the $H$-ideal
generated by $B$ in $R$. Therefore, $(B) \subseteq (0:M)_R$ and $B
\subseteq (0:M)_I,$  which implies $M \in {\cal M}_I$.

Finally, we show that $(M4)$ holds.  Let $I \lhd _H R$ and $\bar R
= R/I$. If $M \in {\cal M}_R$ and $I \subseteq (0:M)_R,$ then $M$
is an $\bar R$-$H$- prime module. If $(N, \bar B)$ has
"$L$-condition"  for $H$-ideal $\bar B$ of $\bar R$ and $\bar
R$-$H$-submodule $N$ of $M$, then subset $F\subseteq B$ and there
exists a natural number $n$ such that $F^nN = (\bar F)^nN =0.$
Consequently, $M\in {\cal M}_{\bar R}$. Conversely, if $M \in
{\cal M}_{\bar R}$, we can similarly show that $M \in {\cal M}_R.$

The second claim follows from Proposition \ref {9.3.18}  and
Theorem \ref {9.3.4}(1). $\Box$

\begin {Theorem}\label {9.3.20}
$r_{Hl}= r_{lH}$.
\end {Theorem}

{\bf Proof.} Obviously, $r_{lH} \le r_{Hl}.$ It remains to show
that $r_{Hl}(R) \not= R$ if $r_{lH}(R) \not=R.$ There exists a
finite subset $F$ of $R$ such that $F^n \not=0$    for  any
natural number $n$. Let $${\cal F} = \{ I \mid I \hbox { is an } H
\hbox {-ideal of } R \hbox { with } F^n \not\subseteq I \hbox {
for any natural number } n \}.$$ By Zorn's lemma, there exists a
maximal element $P$ in ${\cal F}$. It is clear that $P$ is an
$H$-prime ideal of $R.$ Now we show that $r_{lH}(R/P)=0$. If
$0\not=B/P$ is an $H$-ideal of $R/P$, then there exists a natural
number $m$ such that $F^m \subseteq B$. Since $(F^m + P)^n
\not=0+P$ for any natural number $n$, we have that $B/P$ is not
locally nilpotent and $r_{lH}(R/P)=0.$  Consequently, $r_{Hl}(R)
\not=R$. $\Box$

In fact, all of the results hold in braided tensor categories
determined by (co)quasitriangular structure.

\end {document}